\documentclass[leqno,12pt]{article}
\usepackage{amsmath, amssymb}
\usepackage{amsthm}

\theoremstyle{plain}

\theoremstyle{definition}



\newcommand{\beeq}{\begin{equation}}
\newcommand{\eneq}{\end{equation}}

\newcommand{\ba}{\begin{array}}
\newcommand{\ea}{\end{array}}

\newcommand{\be}{\begin{equation}}
\newcommand{\ee}{\end{equation}}
\newcommand{\bea}{\begin{eqnarray}}
\newcommand{\eea}{\end{eqnarray}}

\newcommand{\supp}{\mathrm{supp}}
\date{}
\begin{document}

\baselineskip=9mm

\title{Energy estimates of harmonic maps between Riemannian manifolds}
\author{M. A. Ragusa$^a$$^b$\footnote{Corresponding author: maragusa@dmi.unict.it} and A. Tachikawa$^c$\footnote{tachikawa$ \_ $atsushi@ma.noda.tus.ac.jp}}

\maketitle
\begin{center}
$^a$Dipartimento di Matematica e Informatica, Universit\`{a} di Catania, Viale Andrea Doria,
 6-95125 Catania, Italy\\
$^b$RUDN University, 6 Miklukho - Maklay St, Moscow, 117198, Russia\\
$^c$Department of Mathematics, Faculty of Science and Technology,
Tokyo University of Science, Noda, Chiba, 278-8510, Japan\\
\end{center}

\begin{abstract}
Let  $\Omega \subset {R}^n,$ $n \geq 3,$ be a bounded open set,
$x=(x_1,x_2,\ldots,x_n)$ a generic point which belongs to $\Omega,$
$u \colon \Omega \to {R}^N ,$ $N>1,$
and $ Du=(D_\alpha u^i)$, $D_\alpha = \partial/\partial x_\alpha, $ $\alpha =1,\ldots,n,\,$ $i=1,\ldots,N .\,$

Main goal is the study of regularity of the minima of nondifferentiable functionals
$$
{\cal F} \,=\, \int_\Omega F(x,u,Du) dx.
$$
having the integrand function different shapes of smoothness. The method is based on the use some majorizations for the functional, rather than the well known Euler equation associated to it.
\end{abstract}

\vskip .3cm
Let  $\Omega \subset {R}^n,$ $n \geq 3,$ be a bounded open set,
$x=(x_1,x_2,\ldots,x_n)$ a generic point which belongs to $\Omega,$
$u \colon \Omega \to {R}^N ,$ $N>1,$ $u(x)\,=\, (u^1(x),
\ldots, u^N(x)) $ and $ Du=(D_\alpha u^i)$, $D_\alpha = \partial/\partial x_\alpha, $ $\alpha =1,\ldots,n,\,$ $i=1,\ldots,N .\,$We are interested in regularity of the minima of nondifferentiable functionals
$$
{\cal F} \,=\, \int_\Omega F(x,u,Du) dx.
$$

{\bf Definition}$\,1\,\,$
(see \cite{RT4}).
A minimizer for the functional $ {\cal F}$ is a function $u \in
W^{1,p}(\Omega)$ 
 such that for every
$\varphi \in W^{1,p}_0(\Omega): $
$$
{\cal F}(u;\supp \varphi) \leq  {\cal F}(u+\varphi; \supp \varphi).
$$

Before the remarkable goal achieved by Giaquinta and Giusti in  \cite{GGActa} all regularity results for minima of regular functionals have as starting point the Euler equation of the above functional. This method presents some inconveniences.

At first it requires some smoothness of $F $ and suitable growth conditions, not only on $F, $ but also on its partial derivatives $F_u, $ $F_p. $

Secondly, under natural growth conditions, it is necessary to start with the bounded minimum points $u, $ and also assume, in the vector value case, some smallest condition on $u. $ This, sometimes, does not allows to apply the results to minimum points which are able to find in general only in $H^1_2 . $

From these facts we understand that the Euler method does not distinguish between true minima and simple extremals.

Then, starting with the mentioned paper \cite{GGActa} and the study continued by Giaquinta and Giusti in \cite{GG116 GiaqGIALLO} and \cite{GG117GuaqGIALLO}, was create and built  a regularity theory for minimum points, working directly with the functional ${\cal F} $ instead of acting with its Euler equation, essentially based on a perturbation argument.

We want to underline that later \cite{GG116 GiaqGIALLO} and \cite{GG117GuaqGIALLO} and through the work of many authors, among others  Giusti in \cite{Giusti2003} and Dan\v{e}\v{c}ek and Viszus in \cite{DV}, was developed the study of partial regularity for solutions of nonlinear elliptic systems.

We focus our attention in the case that the coefficients of the principal term are discontinuous. In  \cite{DV} the authors study regularity of minimizers for the functional
$$
\int_\Omega \left\{  A^{\alpha\beta}_{ij}(x) D_\alpha u^i D_\beta u^j + g(x,u,Du)\right\} dx,
$$
being $g$ a lower order term such that for a. e. $x\in \Omega :$
\[
|g(x,u,z)|\, \leq \, f (x) \,+\, K \,|z|^\tau ,
\]
with $f \geq 0,~f \in L^p(\Omega)$,
$ 2\, < \,p \,\leq\   \infty $  $K \geq 0, $ and $0 \leq \tau < 2. $

They prove H\"older regularity of minimizers by assuming that $A^{\alpha\beta}_{ij}(x) $ are in the following vanishing mean oscillation class, at first considered in \cite{S}. Before its definition we need to introduce the space of bounded mean oscillation functions (see \cite{J-N}). We stress that the functions belonging in these classes, used a lot by the authors (see for instance \cite{P-R}, \cite{RT1}, \cite{RT2}, \cite{RT3}), could be discontinuous.

{\bf Definition}$\,2.\,$
  Let $f \in L^1_{\mathrm{loc}}({R}^n).$ We say that $f$ belongs to $BMO({R}^n)$ if is finite the seminorm
\[
\| f \|_* \equiv
\sup_{B(x,R) } \frac{1}{|B(x,R)|}
\int\limits_{B(x,R)} |f(y)-f_{x,R}| dy,
\]
where $B(x,R)$ ranges over the class of the balls of $R^n$  having radius $R > 0$ and center $x \in R^n $.

{\bf Definition}$\,3.\,\,$
Let $f \in BMO({R}^n)$  and
$$
	\eta (f,R) = \sup_{\rho \leq R } \frac{1}{|B(x,\rho)|}
	\int\limits_{B(x,\rho)} |f(y)-f_{\rho}| dy
$$
where  $B(x,\rho)$ is defined as above. We say that $f \in VMO(\Omega)$ if
$$
\lim_{R \to 0} \eta (f,R) =0.
$$
The result contained in \cite{DV} is extended in \cite{RT1}. Indeed in the last mentioned paper the authors treat the functional whose integrand contains the term $g(x,u,Du) $ and has coefficients $A^{\alpha\beta}_{ij} $ dependent not only on $x$ but also on $u. $ Precisely the main goal in \cite{RT1} is the following.

{\bf Theorem 1}.
Let $u \in W^{1,2}(\Omega
)\,$ be a minimum of the functional
$$
{\cal A} \,=\,\int_\Omega \left\{  A^{\alpha\beta}_{ij}(x,u) D_\alpha u^i D_\beta u^j + g(x,u,Du)\right\} dx,
$$
where  $A^{\alpha \beta}_{ij}$ are bounded functions on
$\Omega \times R^N$ and satisfy the conditions
\begin{enumerate}
\item[(A-1)]
     $A^{\alpha\beta}_{ij} = A^{\beta\alpha}_{ji}$.
\item[(A-2)]  For every $u \in R^N$, $A^{\alpha\beta}_{ij}(\cdot , u)
    \in VMO(\Omega)$.
\item[(A-3)]  For every $x\in \Omega$ and $u,v\in R^N$,
    \[
    \big|\,A^{\alpha\beta}_{ij}(x,u)\,-\,A^{\alpha\beta}_{ij}(x,v)\,\big|
       \,\leq\, \omega\,(\,|u-v|^2\,)
    \]
    for some monotone increasing concave function $\omega$ with
    $\lim_{t\to 0+} \omega(t) =0$.
\item[(A-4)] There exist positive constant $\nu$
    such that
    \[
       \nu |\xi|^2 \leq A^{\alpha\beta}_{ij}(x,u)\xi^i_\alpha \xi^j_\beta
    \]
    for a.e. $x\in\Omega$, all $u\in R^N$ and $\xi \in R^{nN}$.
\end{enumerate}

Let us also suppose that the function $g$ is a Charath\'eodory function, that is:
\begin{enumerate}
\item[(g-1)]
$g(\cdot ,u,z )$ is measurable in $x$ $\forall u \in R^N,
\forall z \in R^{nN};$
\item[(g-2)]  $g(x,\cdot , \cdot )$ is continuous in $(u,z)$ $\, a. \,\,e. \, x \in \Omega;$
\end{enumerate}
\begin{enumerate}
\item[(g-3)]
For a. e. $x\in \Omega$
\[
|g(x,u,z)|\, \leq \, g_1 (x) \,+\, H \,|z|^\gamma ,
\]
with $g_1 \geq 0,~g_1 \in L^p(\Omega)$,
$p \in ]2, + \infty[$, $H \geq 0, $ and $0 \leq \gamma < 2.$
\end{enumerate}

Then, if $\lambda = n(1- \frac{2}{p}), $
it follows
\begin{equation}
D\,u \,\in L^{2, \lambda}_{\mathrm{loc}}(\Omega_0
)
\end{equation}
where
\[
\Omega_0\,=\,\{ x \in \Omega \colon \liminf_{R \to 0} \frac
{1}{R^{n-2}}\int_{B(x,R)} |Du(y)|^2 dy \,=\,0
\}.
\]

Moreover, if $n-2<\lambda<n , $ 
 for
$\alpha = 1-\frac{n-\lambda}{2}\in (0,1),$
we have
\begin{equation}
u \in C^{0,\alpha}(\Omega_0
).
\end{equation}

Let us now introduce the notion of Hausdorff measure ( see e. g. \cite{Gia-libro2})

{\bf Definition 4. } Let $X$ be a metric space and $J$ a family of subsets of $X$ such that $ \emptyset\in J. $
Given a function $\zeta : J \to [0, \infty], $ $\zeta (\emptyset) = 0 $ we set  for all subsets $E \subset X$
$$
\mu_\epsilon (E) \,=\, \inf \{\, \sum_{ n = 0}^{\infty} \zeta (F_n)\,:\,\,\, F_n \in J, \,\,\,E \subset \bigcup_{n = 0}^{\infty} F_n, \quad
 \mathrm{diam}
 F_n < \epsilon \,\}
$$
and
$$
\mu (E)\,=\, \lim_{\epsilon \to 0} \mu_\epsilon (E)
$$
where $\mu$ is an exterior measure and is called {\it the result of Charateodory's construction } for $\zeta $ and $J. $ The k-dimensional Hausdorff-measure corrisponds to choosing $X \,=\,R^n, $
$J\,\,$ the open sets in $R^n,  $ $\zeta (F)\,=\, \omega_k 2^{- k} (
\mathrm{diam}
F)^k,  $ where $\omega_k $ is the measure of the unit ball in $R^n,$
$$
{\cal H}^k (E) \,=\, 2^{- k} \omega_k
\sup_{\epsilon \to 0}\big( \inf
\{\,\, \sum_{n = 0}^{\infty}(
\mathrm{diam}
F_n)^k \,\,\}
\big)
$$
being $\{F_k\} $ is a countable family of open sets, $E \subset \bigcup_{n = 0}^{\infty}, $ $\,
\mathrm{diam}
F_n < \epsilon .$

Also we recall that ${\cal H}^0 (E)  $ is the number of points of $E.$

For general quadratic growth functionals, as
$$
{\cal B} \,=\, \int_\Omega B(x,u,Du) dx,
$$
in \cite{RT2} the authors treat the regularity problem of minimizers $u(x) $ assuming that $B(x,u,Du) $
is in the $VMO $ class as function of $x. $
The Morrey regularity result for local minimizers of the functional ${\cal B} $  is contained in the following theorem.

{\bf Theorem 2}.
Let $B(x,u,p)$ be a nonnegative function defined on
$\Omega\times R^N\times R^{mN}$ which satisfies the following conditions.
\begin{enumerate}
\item[(B
-1)]  For every $(u,p)\in R^N \times R^{nN}$,
    $B
    (\cdot, u, p) \in VMO(\Omega)$ and the mean oscillation of
    $B
    (\cdot, u, p)/|p|^2 $ vanishes uniformly with respect to
    $u, p$ in the following sense:
    For some nonnegative function $\sigma(y,\rho)$ with
    \begin{equation}
    \lim_{R\to 0} \sup_{\rho <R}
    \frac{1}{B(x,\rho)}\int
    \sigma(y,\rho) dy =0,
    \end{equation}
    $B
    (\cdot , u ,p)$ satisfies
    \begin{equation}
        \big|B
        (y,u,p) -B
        _{x,\rho}(u,p) \big|
    \leq \sigma(y,\rho)|p|^2 ~~~~ \forall (u,p)
        \in R^N \times R^{nN},
    \end{equation}
    where
    \[
        B
        _{x,\rho}(u,p)=\frac{1}{B(x,\rho)} \int
        B
        (y,u,p) dy.
    \]
\item[(B
-2)]  For every $x\in \Omega$, $p\in R^{nN}$ and $u,v\in R^N$,
    \[ \big|B
    (x,u,p)-B
    (x,v,p)\big|
       \leq \omega(|u-v|^2)|p|^2
    \]
    for some monotone increasing concave function $\omega$ with $\omega(0)
    =0$,
\item[(B
-3)] For almost all $x\in \Omega$ and all $u\in R^N$,
    $B
    (x,u,\cdot) \in C^2(R^{nN})$,
\item[(B
-4)] There exist positive constants $\mu_0 \leq \mu_1$ such that
     \[
     \mu_0|p|^2 \leq B
     (x,u,p) \leq \mu_1 |p|^2
     ~~~~~\mathrm{for~all}~~ (x,u,p)\in \Omega\times R^N \times R^{nN}.
     \]
\end{enumerate}
If $u \in W^{1,2}(\Omega
)\,$ is a local minimizer of the functional
${\cal B
}\, $
 and are true the assumptions \rm (B
-1), (B
-2), (B
-3)
 $\,$and \rm (B
-4) 
then , for every
$0<\lambda <\min\{2+\varepsilon, m\}$ and some $\varepsilon >0, $
we have
\begin{equation}\label{Morrey-reg-u}
D\,u \,\in L^{2, \lambda}_{\mathrm{loc}}(\Omega_0
)
\end{equation}
where
$$
\Omega_0\,=\,\{ x \in \Omega \colon \liminf_{R \to 0} \frac
{1}{R^{n-2}}\int_{B(x,R)} |Du(y)|^2 dy = 0
\}.
$$
Moreover, we have
\[
{\cal H}^{n-2-\delta}(\Omega \setminus \Omega_0)=0
\]
for some $\delta >0$,
where ${\cal H}^r$ is the $r$-dimensional
Hausdorff measure.
$$\,$$
As a consequence of the above theorem we have  that
if $m \leq 4, $
$\,$ for some $\alpha \in (0,1)$, 
\begin{equation}
u \in C^{0,\alpha}(\Omega_0
).
\end{equation}

The investigation has in \cite{RT4} an improvement, the authors consider
$$
{\cal C} \,=\, \int_\Omega C(x,u,Du) dx.
$$
and prove regularity properties in the vector value case for minimizers of variational integrals where the integrand
is not necessarily continuous respect to the variable $x$ but grows polinomially like $| \xi |^p,  $ for $p \geq 2. $ Then this note generalizes \cite{RT2} because, instead of the quadratic growth
assumption, the authors suppose that for some $p \geq 2$  there exist two constants  $\lambda_1$ and $\Lambda_1 $
such that
\begin{equation}
     \lambda_1 (
     1+|\xi|^p) \leq C
     (x,u,\xi) \leq \Lambda_1(
     1 + |\xi|^p), \quad \forall (x,u,\xi)\in \Omega\times R^N \times
     R^{nN},
\end{equation}
Let us now state the mentioned generalization.

{\bf Theorem 3. }
Assume that $\Omega\subset R^n$ is a bounded domain with
sufficiently smooth boundary $\partial \Omega$ and that $p\geq 2$.
Let $u \in H^{1,p}(\Omega
)\,$ be  a minimizer of the functional
\[ {\cal C}
 = \int_\Omega C(x,u,Du) dx
\]
in the class
\[
X_g (\Omega) = \{u \in H^{1,p}(\Omega)~;~
u-g \in H^{1,p}_0(\Omega)\}
\]
for a given boundary data $g \in H^{1,s}(\Omega)$ with $s>p$.
Suppose that the following assumptions 
are satisfied :
\begin{enumerate}
\item[(C
-1)]  For every $(u,\xi)\in R^N \times R^{nN}$,
    $C
    (\cdot, u, \xi) \in VMO(\Omega)$ and the mean oscillation of
    $C
    (\cdot, u, \xi)/|\xi|^p $ vanishes uniformly with respect to
    $u, \xi$ in the following sense:
    there exist a positive number $\rho_0$ and a function
    $\sigma(z, \rho) : R^n \times [0, \rho_0) \to [0,\infty)$ with
    \begin{equation}
    \lim_{R\to 0} \sup_{\rho <R}
    \frac {1}{| Q(0,\rho) |}
    \int
    _{Q(0,\rho)\cap \Omega} \sigma(z,\rho) dz =0,
    \end{equation}
    such that
    $C
    (\cdot , u ,\xi)$ satisfies
    for every $x\in \overline{\Omega}$ and
    $y \in Q(x,\rho_0)\cap \Omega$
    \begin{equation}
	    \big|C
(y,u,\xi) -C
        _{x,\rho}(u,\xi)
        \big| \leq \sigma(x-y,\rho)(
        1 + |\xi|^2)^{\frac{ p}{2}}
          ~~~~ \forall (u,\xi)\in R^N \times R^{nN},
    \end{equation}
    where
    \[
        C
        _{x,\rho}(u,\xi)=  \frac{1}{|Q(x,\rho) |}
        \int
        _{Q(x,\rho)\cap \Omega} C
        (y,u,\xi) dy.
    \]
\item[(C
-2)]  For every $x\in \Omega$, $\xi \in R^{nN}$ and $u,v\in R^N$,
    \[ \big|C
    (x,u,\xi)-C
    (x,v,\xi)\big|
       \leq (
       1 + |\xi|^2 )^{\frac{p}{2}}
       \omega(|u-v|^2)
    \]
  where $\omega$ is some monotone increasing concave function
  with $\omega(0)  =0$,
\item[(C
-3)] For almost all $x\in \Omega$ and all $u\in R^N$,
    $C
    (x,u,\cdot) \in C^2(R^{nN})$,
\item[(C
-4)] There exist positive constants $\lambda_1,$ $\Lambda_1$ such that
     \[
     \lambda_1 (
     1+|\xi|^p) \leq
     C(x,u,\xi) \leq \Lambda_1(
     1 + |\xi|^p)
     \]
     \[
           \lambda_1 (1+|\eta|^p) \leq
           \frac{\partial^2
           C(x,u,\xi)}{\partial \xi^i_\alpha
         \partial \xi^j_\beta} \eta^i_\alpha \eta^j_\beta \leq \Lambda_1(
           1 + |\eta|^p)
         \]
     for all $(x,u,\xi,\eta)\in \Omega\times R^N \times R^{nN}
\times R^{nN}$.
\end{enumerate}
Then, for some positive $\varepsilon$, for every
$0<\tau<\min\{2+\varepsilon, m(1- \frac{p}{s})\}$
we have
\begin{equation}
D\,u \,\in L^{p, \tau}(\Omega_0
)
\end{equation}
where $\Omega_0$ is a relatively open subset of $\overline{\Omega}$ which
satisfies
$$
\overline{\Omega}\setminus\Omega_0
= \{ x \in \Omega \colon \liminf_{R \to 0} \frac
{1}{R^{n-p}}\int_{\Omega(x,R)} |Du(y)|^p dy > 0
\}.
$$
Moreover, we have
\[
{\cal H}^{n- p -\delta}(\overline{\Omega} \setminus \Omega_0)=0
\]
for some $\delta >0. $

As a corollary of the above theorem the following partial
H\"{o}lder regularity result is true.

{\bf Corollary 1}
Let
$g$, $u$ and $\Omega_0$ be as in Theorem 3. 
 Assume that $p+2\geq n
 $
and that $s > \max\{n
,p\}$.
Then, for some $\alpha \in (0,1)$, we have
\begin{equation}\label{3.2}
u \in C^{0,\alpha}(\Omega_0
).
\end{equation}
Moreover, as a consequence 
 it follows the following full-regularity result if  the integrand $C
$ does not depend on $u$.

{\bf Corollary 2}
Assume that $C
$ and $g$ satisfy all assumptions of Theorem 3 
and that the integrand of the functional ${\cal C}$ 
 does not depend on $u$.
Let $u$ its minimizer 
 in the class $X_g (\Omega), $
then
\begin{equation}\label{Morrey-reg-u-2}
D\,u \,\in L^{p, \tau}(\Omega
).
\end{equation}

Moreover, if $p+2\geq n
$ and $s> \max\{n
,p\}, $
 we have full-H\"{o}lder regularity of $u$.
In conclusion we obtain
\[
u \in C^{0,\alpha} (\overline{\Omega}
).
\]
Let us now consider the variational integral
$$
{\cal D} \,=\, \int_\Omega D(x,u,Du) dx.
$$
assuming that $D(x,u,Du)\in \Omega \times R^N\times R^{n N} $ has $p-$growth :
\begin{equation}
	\lambda_0 (\mu^2 + |\xi|^2)^{p/2}
	\leq D(x,u,\xi) \leq \Lambda_0 (\mu^2+|\xi|^2)^{p/2}
\end{equation}
for all $(x,u,\xi)\in \Omega \times R^N \times R^{
nN}$.
Let us also suppose the following convexity condition
\begin{equation}
	\lambda_0 (\mu^2 + |\xi|^2)^{p/2} |\eta|^2
	\leq \frac{\partial^2
	D(x,u,\xi)}{\partial \xi^j_\beta \partial \xi^i_\alpha}
	\eta^i_\alpha \eta^j_\beta \leq \Lambda_1 (\mu^2+ |\xi|^2)^{p/2}
	|\eta|^2
\end{equation}
for all $(x,u,\xi,\eta) \in \Omega \times R^N
 \times R^{
 n N} \times R^{
 n N }$.

If $\mu =0, $ the variational peculiarity of the functional is
dissimilar from the case $\mu \neq 0. $

In \cite{RT4} are proved partial regularity results for the case
$\mu \neq 0, $  open problem that the authors wish to investigate is the case $\mu = 0. $

\bibliographystyle{plain}

\end{document}